\documentclass[11pt]{article}  
 \oddsidemargin=\evensidemargin
 \usepackage{color}              
 \usepackage{epsfig}
 \usepackage{graphicx}
 \usepackage{amssymb}
 \usepackage{epsfig}

 \newcommand{\cd}{\ \stackrel{d}{\rightarrow} \ }

 \newcommand{\BB}{\mathcal{B}}

 \newcommand{\sfrac}[2]{{\textstyle\frac{#1}{#2}}}

 \newcommand{\ego}{{\bf ego}}
 \newcommand{\payoff}{\mathrm{payoff}}
 \newcommand{\reward}{\mathrm{reward}}
 \newcommand{\thnear}{\theta_{\mbox{{\tiny near}}}}
 \newcommand{\thfar}{\theta_{\mbox{{\tiny far}}}}
 \newcommand{\thNnear}{\theta_{N,\mbox{{\tiny near}}}}
 \newcommand{\thNfar}{\theta_{N,\mbox{{\tiny far}}}}

 \newcommand{\phiNnear}{\phi_{N,\mbox{{\tiny near}}}}
 \newcommand{\phiNfar}{\phi_{N,\mbox{{\tiny far}}}}
 
\begin{document}
 \author{David J. Aldous\thanks{Research supported by N.S.F. Grant
 DMS0704159}
 \\
 \\
        University of California\\
        Department of Statistics\\
         367 Evans Hall \# 3860\\
        Berkeley CA 94720-3860}

 \title{When Knowing Early Matters: Gossip, Percolation and Nash
 Equilibria}

 \maketitle

 \begin{abstract}
 Continually arriving information is communicated through a network of
 $n$
 agents, with the value of information to the $j$'th recipient being a
 decreasing function of $j/n$, and communication costs paid by recipient.
 Regardless of details of network and communication costs, the social
 optimum policy is to communicate arbitrarily slowly.  But selfish agent
 behavior leads to Nash equilibria which
 (in the $n \to \infty$ limit) may be
 efficient (Nash payoff $=$ social optimum payoff)
 or wasteful ($0 < $  Nash payoff $<$ social optimum payoff)
 or totally wasteful (Nash payoff $=0$).
 We study the cases of the complete network (constant communication costs
 between all agents), the grid
 with only nearest-neighbor communication, and the grid with
 communication cost
 a function of distance.
 The main technical tool is analysis of the associated first passage
 percolation process or SI epidemic (representing spread of one item of information) and
 in particular its ``window width", the time interval during which most
 agents learn the item.

 Many arguments are just outlined, not intended as complete rigorous proofs. 
 This version was written in July 2007 to accompany a talk at the ICTP workshop 
 ``Common Concepts in Statistical Physics and Computer Science", and intended as a starting point 
for future thesis projects which could explore  these and many variant problems in detail.  
One of the topics herein 
(first passage percolation on the $N \times N$ torus with short and long range interactions; section 
\ref{sec-5-fpp})
has now been studied rigorously  by  Chatterjee and Durrett \cite{chat-durr}, 
and so it seems appropriate to make this version publicly accessible.

 \end{abstract}

 \newpage
 \section{Introduction}
 A topic which one might loosely call
 ``random percolation of information through networks"
arises in many different contexts, from
  epidemic models \cite{MR1784822} and computer virus models \cite{computer-virus}
  to {\em gossip algorithms}
 \cite{1039491} designed to keep nodes of a decentralized network updated
 about information needed to maintain the network.
This topic differs from {\em communication networks} in that we envisage information as having a definite source but no definite destination.

 In this paper we study an aspect where the vertices of the network
 are agents, and where there are costs and benefits associated with the
 different choices that agents may make in communicating information.  In
 such ``economic game theory" settings one anticipates a {\em social
 optimum} strategy that maximizes the total net payoff to all agents
 combined, and an (often different)
 {\em Nash equilibrium} characterized by the property that no one agent
 can
 benefit from deviating from the Nash equilibrium strategy followed by
 all
 other agents (so one anticipates that any reasonable process of agents adjusting strategies in a
 selfish way will lead to some Nash equilibrium).
 Of course a huge number of different models of costs, benefits and
 choices
 could fit the description above, but we focus on the specific setting
 where
 the value to you of receiving information depends on how few people know the
 information before you do.
 Two familiar real world examples are {\em gossip} in social networks and
 {\em insider trading} in financial markets.
 In the first, the gossiper gains perceived social status from
 transmitting
 information, and so is implicitly willing to pay for communicate to
 others; in the second the owner of knowledge recognizes its value and
 implictly
 expects to be paid for communication onwards.   Our basic model makes
 the
 simpler assumption that the value to an agent attaches at the time
 information is received, and subsequently the agent takes no initiative to
 communicate it to others,
 but does so freely when requested, with the requester paying the cost of
 communication.
 In our model the benefits come from, and communication costs are paid
 to,
 the outside world: there are no payments between agents.

 \medskip
 \noindent
 {\bf Remark.} 
  Many arguments are just outlined, not intended as complete rigorous proofs. 
 This version was written in July 2007, and intended as a starting point 
for future thesis projects which could explore  these and many variant problems in detail.  
One of the topics herein 
(first passage percolation on the $N \times N$ torus with short and long range interactions; section 
\ref{sec-5-fpp})
has now been studied rigorously  by  Chatterjee and Durrett \cite{chat-durr}, 
and so it seems appropriate to make this version publicly accessible.

 \subsection{The general framework: a rank-based reward game}
 \label{sec-frame}
 There are $n$ agents (our results are in the $n \to \infty$ limit). The
 basic two rules are:

\noindent
 {\bf  Rule 1.}
 New items of information arrive at times of a rate-$1$ Poisson process;
 each item comes to one random agent.

 Information spreads between agents by virtue of one agent calling another and learning all items that the other knows (details are case-specific, described later),  with a (case-specific) 
 communication cost paid by the {\em receiver} of information.

\noindent
 {\bf Rule 2.} The $j$'th person to learn an item of information gets reward
 $R(\frac{j}{n})$.

 Here $R(u), \ 0 < u \leq 1$ is a  function such that
 \begin{equation}
  R(u) \mbox{ is decreasing; } R(1) = 0; \quad
 0 < \bar{R}:= \int_0^1 R(u) du < \infty .
 \label{R-props}
 \end{equation}
 Assuming information eventually reaches each agent, the total reward
 from
 each item will be
 $\sum_{j=1}^n R(\frac{j}{n}) \sim n \bar{R}$.
 If agents behave in some ``exchangeable" way then the average net payoff
 (per agent per unit time) is
 \begin{equation}
  \payoff  = \bar{R} - \mbox{(average communication cost per agent per unit
 time)} .
  \label{payoff-R}
  \end{equation}
 Now the average communication cost per unit time can be made arbitrarily
 small
 by simply communicating less often
 (because an agent learns {\em all} items
 that another agent knows, for the cost of one call.  Note the calling agent does not know in advance whether the other agent has any new items of information).
 Thus the ``social optimum" protocol is to communicate arbitrarily slowly, giving
 payoff arbitrarily close to $\bar{R}$.
 But if agents behave selfishly then one agent may gain an advantage by
 paying to obtain information  more quickly, and so we seek to study Nash
 equilibria for selfish agents.
 In particular there are three qualitative different possibilities. In the
 $n \to \infty$ limit, the Nash equilibrium may be 
 
 \noindent
$\bullet$  efficient (Nash payoff $=$ social optimum payoff)\\
 $\bullet$   or wasteful ($0 < $  Nash payoff $<$ social optimum payoff)\\
$\bullet$    or totally wasteful (Nash payoff $=0$).

 \subsection{Methodology}
 Allowing agents' behaviors to be completely general makes the problems
 rather complicated
 (e.g. a subset of agents could seek to coordinate their actions) so in
 each specific model we restrict
 agent behavior to be of a specified form, making calls at random times with a rate parameter $\theta$; the
 agent's ``strategy" is just a choice of $\theta$, and for this
 discussion
 we assume $\theta$ is  a single  real number.
 If all agents use the same parameter value $\theta$ then the spread of one item of 
 information through the network is as some model-dependent  first
 passage
 percolation process 
(see section \ref{FPP-setup}).
 So there is some function
 $F_{\theta,n}(t)$ giving the proportion of agents who learn the item
 within time $t$ after the arrival of the information into the network. Now
 suppose one agent \ego\ uses a different parameter value $\phi$ and gets
 some payoff-per-unit-time, denoted by
 $\payoff(\phi,\theta)$.
 The Nash equilibrium value   $\theta^{\mbox{{\tiny Nash}}}$ is the value
 of $\theta$ for which \ego\ cannot
 do better by choosing a different value of $\phi$, and hence
 is the solution of
 \begin{equation}
 \left. \frac{d}{d \phi} \payoff(\phi,\theta) \right|_{\phi = \theta} = 0
 .
 \label{Nash-crit}
 \end{equation}
 Obtaining a formula for $\payoff(\phi,\theta)$ requires knowing
 $F_{\theta,n}(t)$
 and knowing something about the geometry of the sets of informed agents
 at
 time $t$ --  see (\ref{payoff-CG},\ref{payoff-NN}) for the two basic examples. The
 important point is that where we know the exact $n \to \infty$
 limit
 behavior of $F_{\theta,n}(t)$ we get a formula for the exact limit
 $\theta^{\mbox{{\tiny Nash}}}$,
 and where we know order of magnitude behavior of $F_{\theta,n}(t)$
  we get order of magnitude behavior of $\theta^{\mbox{{\tiny Nash}}}$.

  Note that we have assumed that  in a Nash equilibrium each agent uses
 the same strategy.
  This is only a sensible assumption when the network cost structure has
 enough symmetry
  (is {\em transitive} -- see section \ref{sec-transitive}) and 
 the non-transitive case 
 is an interesting topic for future study.

 It turns out (section \ref{sec-heuristics}) that for determining the qualitative behavior of the Nash
 equilibria, the important aspect is the size of the
 {\em window width} $w_{\theta,n}$ of the associated
  first passage percolation process,
 that is the time interval over which the proportion of agents knowing the
 item of information increases from (say) 10\% to 90\%.
 While this is well understood in the simplest examples of first passage
 percolation on finite sets, it has not been studied for very general
 models and our game-theoretic questions provide motivation for future such
 study.

 To interpret later formulas it turns out to be convenient to work with
 the
 derivative of $R$.
 Write $R^\prime(u) = - r(u)$, so that
 $R(u) = \int_u^1 r(s) ds$  and (\ref{R-props}) becomes
 \begin{equation}
 r(u) \geq 0; \quad
 0 < \bar{R}:= \int_0^1 u r(u) du < \infty .
 \label{r-props}
 \end{equation}

 \subsection{Summary of results}

 \subsubsection{The complete graph case}

 {\bf Network communication  model:}
 Each agent $i$ may, at any time, call any other agent $j$ (at cost $1$),
 and learn all items that $j$ knows.

\noindent
 {\bf Poisson strategy.}
 The allowed strategy for an agent $i$ is to place calls, at the times of
 a
 Poisson
 (rate $\theta$) process, to a random agent.

\noindent
 {\bf Result} (section \ref{sec-cg}).  In the $n \to \infty$ limit the
 Nash equilibrium value of
 $\theta$ is
 \begin{equation}
  \theta^{\mbox{{\tiny Nash}}} =
 \int_0^1 (1 + \log(1-u)) R(u) du
 = \int_0^1 r(u) g(u) du
  \label{CG-result}
 \end{equation}
 where
 $g(u) = -(1-u) \log(1-u) > 0$.

 Our assumptions (\ref{R-props}) on $R(u)$ imply
 $0 < \theta^{\mbox{{\tiny Nash}}} < \bar{R}$.
 Because an agent's average cost per unit time equals his value of
 $\theta$,
 from (\ref{payoff-R})
 the Nash equilibrium payoff $\bar{R} - \theta^{\mbox{{\tiny Nash}}}$ is
 strictly less than the social optimum payoff
  $\bar{R} $ but strictly greater than $0$.
  So this is a ``wasteful" case.

  \subsubsection{The nearest neighbor grid}
  {\bf Network communication  model:}
 Agents are at the vertices of the $N \times N$ torus
 (i.e. the grid with periodic boundary conditions).
 Each agent $i$ may, at any time, call any of the $4$ neighboring agents
 $j$ (at cost $1$), and learn all items that $j$ knows.

\noindent
 {\bf Poisson strategy.}
 The allowed strategy for an agent $i$ is to place calls, at the times of
 a
 Poisson
 (rate $\theta$) process, to a random neighboring agent.

\noindent
 {\bf Result} (section \ref{sec-NN}).
 The Nash equilibrium value of $\theta$ is such that as $N \to \infty$ 
 \begin{equation}
  \theta^{\mbox{{\tiny Nash}}}_N \sim N^{-1}
 \int_0^1 g(u) r(u) du
  \label{torus-result}
 \end{equation}
 where $g(u) > 0$ is a certain complicated function -- see (\ref{gz-formula}).

 So here the Nash equilibrium payoff $\bar{R} - \theta^{\mbox{{\tiny
 Nash}}}_N$
 tends to $\bar{R}$; this is an ``efficient" case.

 \subsubsection{Grid with communication costs increasing with distance}
 \label{sec-sec-cc}
 {\bf Network communication  model.}
  The agents are at the vertices of the $N \times N$ torus.
  Each agent $i$ may, at any time, call any other agent $j$,
  at cost $c(N,d(i,j))$,
  and learn all items that $j$ knows.

  Here $d(i,j)$ is the distance between $i$ and $j$.
 We treat two cases, with different choices of $c(N,d)$.
 In section \ref{sec-grid-cN} we take cost function
 $c(N,d) = c(d)$ satisfying
 \begin{equation}
 c(1) = 1; \quad c(d) \uparrow \infty \mbox{ as } d \to \infty
 \label{cdcd}
 \end{equation} 
 and
 
 \noindent
 {\bf Poisson strategy.}
 An agent's strategy is described by a sequence
 $(\theta(d); \ d = 1,2,3,\ldots)$;  where for each $d$:\\
 \hspace*{0.7in}
 at rate $\theta(d)$ the agent calls a random agent at distance $d$.

 In this case a simple abstract argument (section \ref{sec-grid-cN}) shows
 that the Nash equilibrium is efficient (without calculating what the 
 equilibrium strategy or payoff actually is) for any $c(d)$ satisfying  (\ref{cdcd}).

 In section \ref{sec-torus-general} we take
 \begin{eqnarray*}
 c(N,d) &=& 1; \quad d = 1\\
 &=& c_N; \quad d > 1
 \end{eqnarray*}
 where $1 \ll c_N \ll N^3$, and

\noindent
 {\bf Poisson strategy.}
 An agent's strategy is described by a pair of numbers
 $(\thnear, \thfar) = \theta$: \\
 \hspace*{0.7in}
 at rate $\thnear$ the agent calls a random neighbor \\
 \hspace*{0.7in}
 at rate $\thfar$ the agent calls a random non-neighbor.

 In this case we show (\ref{eq-99}) that the Nash equilibrium strategy satisfies
\[
 \theta_{\mbox{{\tiny near}}}^{\mbox{{\tiny Nash}}} \sim \zeta_1
 c_N^{-1/2}; \quad
 \theta_{\mbox{{\tiny far}}}^{\mbox{{\tiny Nash}}} \sim \zeta_2
 c_N^{-2} \]
 for certain constants $\zeta_1, \zeta_2$ depending on the reward function. 
So the Nash equilibrium cost $\sim \zeta_1 c_N^{-1/2}$,
 implying that the equilibrium is efficient.

 \subsubsection{Plan of paper}
 The two basic cases (complete graph, nearest-neighbor grid) can be
 analyzed directly using known results for first passage percolation on
 these structures; we do this analysis in sections \ref{sec-cg} and
 \ref{sec-NN}.
 There are of course simple arguments for order-of-magnitude behavior in
 those cases, which we recall in section \ref{sec-heuristics} (but which
 the reader may prefer to consult first) as a preliminary to the more
 complicated model ``grid with communication costs increasing with distance", for which one needs to understand orders of magnitude
 before embarking on calculations.

 \subsection{Variant models and questions}
 These results suggest many alternate questions and models,
 a few of which are addressed briefly in the sections indicated, the others
 providing suggestions for future research.

 \begin{itemize}
 \item Are there cases where the Nash equilibrium is totally wasteful?
 (section \ref{sec-finite})
 \item Wouldn't it be better to place calls at regular time intervals?
 (section \ref{sec-reg})
\item Can one analyze more general strategies?
 \item In the grid context of section \ref{sec-sec-cc}, what is the
 equilibrium strategy and cost for more general costs $c(N,d)$?
 \item  What about the symmetric model where, when $i$ calls $j$, they
 exchange information? (section \ref{sec-transitive})
 \item In formulas (\ref{CG-result},\ref{torus-result})  we see
 decoupling between the reward function $r(u)$ and the function $g(u)$
 involving the rest of the model -- is this a general phenomenon?
 \item In the nearest-neighbor grid case, wouldn't it be better to cycle
 calls through the $4$ neighbors?
 \item What about non-transitive models, e.g. social networks where
 different agents have different numbers of friends, so that different
 agents have different strategies in the Nash equilibrium?
 \item To model gossip, wouldn't it be better to make the reward to agent
 $i$ depend on the number of other agents who learn the item from agent
 $i$? (section \ref{sec-audience})
 \item To model insider trading, wouldn't it be better to say that agent
 $j$ is willing to pay some amount
 $s(t)$ to agent $i$ for information that $i$ has possessed for time $t$, the
 function $s(\cdot)$ not specified in advance but a component of strategy
 and hence with a Nash equilibrium value?
 \end{itemize}

 \subsection{Conclusions}
 As the list above suggests, we are only scratching the surface of a
 potentially large topic.  In the usual setting of information
 communication networks, the goal is to communicate quickly, and our two
 basic examples
 (complete graph; nearest-neighbor grid) are the extremes of rapid and slow
communication.  It is therefore paradoxical that, in our rank-based reward
 game, the latter is efficient while the former is inefficient.
 One might jump to the conclusion that in general efficiency in the
 rank-based reward game was inversely related to network connectivity.  But
 the examples of the grid with long-range interaction show the
 situation is not so simple, in that agents {\em could} choose to make long
 range calls and emulate a highly-connected network, but in equilibrium
 they do not do so  very often.

 \section{The complete graph}
 \label{sec-cg}
 The default assumptions in this section are

\noindent
 {\bf Network communication  model:}
 Each agent $i$ may, at any time, call any other agent $j$ (at cost $1$),
 and learn all items that $j$ knows.

\noindent
 {\bf Poisson strategy.}
 The allowed strategy for an agent $i$ is to place calls, at the times of
 a
 Poisson
 (rate $\theta$) process, to a random agent.

 \subsection{Finite number of rewards}
 \label{sec-finite}
 Before deriving the result (\ref{CG-result}) in our general framework,
 let
 us step outside that framework
 to derive a very easy variant result.
 Suppose that only the first two recipients of an item of information receive
 a
 reward, of amount $w_n$ say.
 Agent strategy cannot affect the first recipient, only the second.
 Suppose \ego\ uses rate $\phi$ and other agents use rate $\theta$. Then
 (by elementary properties of Exponential distributions)
 \begin{equation}
 P( \mbox{\ego\  is second to receive item}) = \frac{\phi}{\phi +
 (n-2)\theta}
 \label{k=2ego}
 \end{equation}
 and so
 \[ \payoff(\phi,\theta) = \frac{w_n}{n} + \frac{\phi w_n}{\phi +
 (n-2)\theta} - \phi . \]
 We calculate
 \[ \frac{d}{d\phi} \payoff(\phi,\theta) = \frac{(n-2)\theta w_n}{(\phi +
 (n-2)\theta)^2} - 1 \]
 and then the criterion (\ref{Nash-crit}) gives
 \[ \theta^{\mbox{{\tiny Nash}}}_n = \frac{(n-2)w_n}{(n-1)^2} \sim
 \frac{w_n}{n} . \]
 To compare this variant with the general framework, we want the total
 reward available from an item to equal $n$, to make  the social optimum
 payoff $\to 1$, so we choose $w_n = n/2$.
 So we have shown that the Nash equilibrium payoff is
 \begin{equation}
  \payoff = 1 - \theta^{\mbox{{\tiny Nash}}}_n \to \sfrac{1}{2} .
   \label{nash-2-1}
  \end{equation}
 So this is a ``wasteful" case.

 By the same argument
 we can study the case where (for fixed $k \geq 2$) the first $k$
 recipients get reward
 $n/k$.
 In this case we find
 \[ \theta^{\mbox{{\tiny Nash}}}_n  \sim \frac{k-1}{k}  \]
 and the Nash equilibrium payoff is
 \begin{equation}
  \payoff \to \sfrac{1}{k}
  \label{nash-2-2}
  \end{equation}
 while the social optimum payoff $= 1$.
 Thus by taking $k_n \to \infty$ slowly we have a model in which the Nash
 equilibrium is
 ``totally wasteful".

 \subsection{First passage percolation : general setup}
\label{FPP-setup}
The classical setting for first passage 
 percolation, surveyed in \cite{MR2045986},
concerns nearest neighbor percolation on the $d$-dimensional lattice.
 Let us briefly state our general setup 
 for first passage percolation (of ``information")
 on a finite graph.  
 There are ``rate" parameters $\nu_{ij} \geq 0$ for undirected edges $(i,j)$. 
 There is an initial vertex $v_0$, which receives the information at time $0$.  
 At time $t$, for each vertex $i$ which has already received the information, and each neighbor $j$, 
 there is chance $\nu_{ij} dt$ that $j$ learns the information from $i$ before time $t+ dt$.  
 Equivalently, create independent Exponential($\nu_{ij}$) random variables $V_{ij}$ on edges $(i,j)$.
 Then each vertex $v$ receives the information at time
 \[ T_v = \min \{V_{i_0i_1} + V_{i_1i_2} +  \ldots + V_{i_{k-1}i_k} \} \]
 minimized over paths $v_0 = i_0, i_1, i_2, \ldots, i_k = v$.

 \subsection{First passage percolation on the complete graph}
 Consider first passage
 percolation on the complete $n$-vertex graph with rates $\nu_{ij} = 1/(n-1)$.  
 Pick $k$ random agents and
 write $\bar{S}^n_{(1)},\ldots,\bar{S}^n_{(k)}$ for the times
 at which these  $k$  agents receive the information.  The key fact for our
 purposes is that as $n \to \infty$
 \begin{equation}
 (\bar{S}^n_{(1)} - \log n,\ldots,\bar{S}^n_{(k)} - \log n) \cd
 (\xi + S_{(1)}, \ldots, \xi + S_{(k)})
 \label{KS-lim}
 \end{equation}
 where the limit variables are independent, $\xi$ has double exponential
 distribution
 $P(\xi \leq x) = \exp(-e^{-x})$
 and each $S_{(i)}$ has
 the {\em logistic} distribution with distribution
 function
 \begin{equation}
  F_1(x) = \frac{e^x}{1+ e^x} , \quad - \infty < x < \infty .
 \label{def-logistic}
 \end{equation}
 Here $\cd$ denotes convergence in distribution.
 To outline a derivation of (\ref{KS-lim}), fix a large integer $L$  and
 decompose the percolation times as
 \begin{equation}
  \bar{S}^n_{(i)} - \log n
 = (\tau_L - \log L) +
  (\bar{S}^n_{(i)} - \tau_L + \log (L/n))  
  \label{n-decomp}
  \end{equation}
 where $\tau_L$ is the time at which some $L$ agents have received the
 information. 
 By the Yule process approximation
 (see e.g. \cite{me-shanky})
 to the fixed-time behavior of the first passage percolation,
 the number $N(t)$ of agents possessing the information at fixed large time $t$
 is approximately distributed as $We^t$, where $W$ has Exponential($1$) distribution, and so
 \[ P(\tau_L \leq t)
 = P(N(t) \geq L)
 \approx P(We^t \geq L)
 = \exp(-Le^{-t}) \]
 implying
 $\tau_L - \log L \approx \xi$ in distribution,
 explaining the first summand on the right side of (\ref{KS-lim}).   Now
 consider the proportion $H(t)$ of agents possessing the information at
  time $\tau_L + t$.
 This proportion follows closely the deterministic logistic equation
 $H^\prime = H(1-H)$ whose solution is (\ref{def-logistic}) shifted to
 satisfy the initial condition $H(0) = L/n$, so this solution approximates
 the distribution function of $S_{(i)} - \log (L/n)$.
 Thus the time $\bar{S}^n_{(i)}$ at which a random agent receives the
 information satisfies
 \[  (\bar{S}^n_{(i)} - \tau_L + \log (L/n)) \approx S_{(i)}
 \mbox{ in distribution } \]
 independently as $i$ varies.  Now the limit decomposition (\ref{KS-lim}) follow from the finite-$n$ 
 decomposition (\ref{n-decomp})..

 We emphasize (\ref{KS-lim}) instead of  more elementary derivations
 (using methods of \cite{MR0100305,MR1965974})
 of the limit distribution for $\bar{S}^n_{(1)} - \log n$
 because (\ref{KS-lim}) gives the correct dependence structure for
 different agents.  Because only {\em relative} order of gaining
 information is relevant to us, we may recenter by subtracting $\xi$ and
 suppose that the times at which different random agents gain information
 are independent
 with logistic distribution (\ref{def-logistic}).

 \subsection{Analysis of the rank-based reward game}
 \label{sec-CG-analy}
 We now return to our general reward framework
 \begin{quote}
 The $j$'th person to learn an item of information gets reward
 $R(\frac{j}{n})$
 \end{quote}
 and give the argument for (\ref{CG-result}).

 Suppose all agents use the Poisson($\theta$) strategy.  
 In the case $\theta = 1$, the way that a single item of information spreads is exactly as the first passage percolation process above; and the general-$\theta$ case is just a time-scaling by $\theta$.
 So as above, we may suppose that
 (all calculations in the $n \to \infty$ limit) the recentered time $S_\theta$ to
 reach a random agent
 has distribution function
 \begin{equation}
  F_\theta(x) = F_1(\theta x) 
  \label{F-theta}
  \end{equation}
 which is the solution of the time-scaled logistic equation
 \begin{equation}
 \frac{F^\prime_\theta}{1 - F_\theta} = \theta  F_\theta 
 \label{eq-logistic-theta}
 \end{equation}
 (Recall $F_1$ is the logistic distribution (\ref{def-logistic})).
 Now consider the case where all other agents use a value $\theta$ but
 \ego\ uses a different value $\phi$.
 The (limit, recentered) time $T_{\phi,\theta}$ at which \ego\ learns the
 information now has
 distribution function
 $G_{\phi,\theta}$ satisfying an analog of (\ref{eq-logistic-theta}):
\begin{equation}
 \frac{G^\prime_{\phi,\theta}}{1 - G_{\phi,\theta}} = \phi  F_\theta . 
 \label{Gphit}
 \end{equation}
 To explain this equation, the left side is the rate at time $t$ at which \ego\  learns the information; this equals the rate $\phi$ of calls by \ego, times the probability $F_\theta(t)$ that the called agent has received the information.
To solve the equation, first we get
 \[ 1 - G_{\phi,\theta} = \exp \left( - \phi \int F_\theta \right) . \]
 But we know that in the case $\phi = \theta$ the solution is $F_\theta$,
 that is we know
 \[ 1 - F_{\theta} = \exp \left( - \theta \int F_\theta \right) , \] and
 so we have the solution of (\ref{Gphit})  in the form
 \begin{equation}
 1 - G_{\phi,\theta} =
 (1 - F_\theta)^{\phi/\theta} .
 \label{eq-G}
 \end{equation}
 If \ego\ gets the information at time $t$ then his percentile rank is
 $F_\theta(t)$ and his reward is $R(F_\theta(t))$.
 So the expected reward to \ego\ is
 \[ E R(F_\theta(T_{\phi,\theta})); \quad
 \mbox{where dist}(T_{\phi,\theta}) = G_{\phi,\theta} . \]
 We calculate
 \begin{eqnarray}
 P(F_\theta(T_{\phi,\theta}) \leq u)
 &=& G_{\phi,\theta}(F_\theta^{-1}(u))\nonumber \\
 &=& 1 - (1 - F_\theta(F_\theta^{-1}(u)) )^{\phi/\theta} \mbox{ by
 (\ref{eq-G})}  \nonumber \\
 &=&
 1 - (1-u)^{\phi/\theta} \label{FTu}
 \end{eqnarray}
 and so 
 \[ E R(F_\theta(T_{\phi,\theta})) =  \int_0^1 r(u) \
 (1 - (1-u)^{\phi/\theta} )du . \]
 This is the mean reward to \ego\ from one item, and hence also the mean reward per unit time in the ongoing process.
So, including the ``communication cost" of $\phi$ per unit time, 
 the net payoff (per unit time) to \ego\ is
 \begin{equation}
 \payoff(\phi,\theta) =
 - \phi
 + \int_0^1 r(u) \
 (1 - (1-u)^{\phi/\theta} )du .
 \label{payoff-CG}
 \end{equation}
 The criterion (\ref{Nash-crit}) for $\theta$ to be a
 Nash equilibrium is, using the fact
 $\frac{d}{d \phi} x^{\phi/\theta} = \frac{\log x}{\theta}
 x^{\phi/\theta}$,
 \begin{equation}
  1 = \sfrac{1}{\theta}
 \int_0^1  r(u)  \
 (- \log (1-u)) \
 (1-u) du .
 \label{eq-u}
 \end{equation}
 This is the second equality in (\ref{CG-result}),
 and integrating by parts gives the first equality.

 \paragraph{ Remark.}
 For the linear reward function
 \[ R(u) = 2(1-u); \quad \bar{R} = 1 \]
 result (\ref{CG-result}) gives Nash payoff $= 1/2$.
 Consider alternatively
 \[ R(u) = \sfrac{1}{u_0} 1_{(u \leq u_0)}; \quad
 \bar{R} = 1 . \]
 Then the $n \to \infty$ Nash equilibrium cost is
 \[  \theta^{\mbox{{\tiny Nash}}}(u_0) =  \frac{1}{u_0}
 \int_0^{u_0} (1 + \log(1-u)) \ du .\]
In particular,  the Nash payoff $1 - \theta^{\mbox{{\tiny Nash}}}(u_0) $ satisfies
 \[ 1 - \theta^{\mbox{{\tiny Nash}}}(u_0) \to 0 \mbox{ as } u_0 \to 0 .
 \]
 In words, as the reward becomes concentrated on a smaller and smaller
 proportion of the population then
 the Nash equilibrium becomes more and more wasteful.
 In this sense
 result (\ref{CG-result}) in the general framework is consistent with the
 ``finite number of rewards" result (\ref{nash-2-2}).

 \section{The $N \times N$ torus, nearest neighbor case}
 \label{sec-NN}

  {\bf Network communication  model.}
  There are $N^2$  agents at the vertices of the $N \times N$ torus.
 Each agent $i$ may, at any time, call any of the $4$ neighboring agents
 $j$ (at cost $1$), and learn all items that $j$ knows.

\noindent
 {\bf Poisson strategy.}
 The allowed strategy for an agent $i$ is to place calls, at the times of
 a
 Poisson
 (rate $\theta$) process, to a random neighboring agent.

 We will derive formula (\ref{torus-result}).
 As remarked later, the function $g(u)$ is ultimately derived from fine
 structure
 of first passage percolation in the plane, and seems impossible to
 determine as an explicit formula.
 But of course the main point is that
 (in contrast to the complete graph case)
 the Nash equilibrium payoff $\bar{R} - \theta^{\mbox{{\tiny Nash}}}_N =
 \bar{R} - O(N^{-1})$
 tends to the social optimum $\bar{R}$.

\subsection{Nearest-neighbor first passage percolation on the torus}
 Consider (nearest-neighbor)  first passage percolation on the $N \times N$ torus, started at
 a
 uniform random vertex, with
rates $\nu_{ij} = 1$ for  edges $(i,j)$.
 Write
 $(T^N_i, \ 1 \leq i \leq 4)$
 for the information receipt times of the $4$ neighbors of the origin (using paths
 not
 through the origin), and write
 $Q^N(t)$ for the number of vertices informed by time $t$.
 Write
 $T^N_* = \min (T^N_i, \ 1 \leq i \leq 4)$.

 The key point is that we expect a $N \to \infty$ limit of the following
 form
 \[
  (T^N_i - T^N_*, 1 \leq i \leq 4; \ N^{-2} Q^N(T^N_*); \
 (N^{-1}(Q^N(T^N_*
 + t) - Q^N(T^N_*)), 0 \leq t < \infty))
 \]
 \begin{equation}
 \cd
 (\tau_i , 1 \leq i \leq 4;\  U; \ (V t, 0 \leq t < \infty))
 \label{torus-con}
 \end{equation}
 where $\tau_i , 1 \leq i \leq 4$ are nonnegative with $\min_i \tau_i =
 0$;
 $U$ has uniform$(0,1)$ distribution;
 $0<V<\infty$;  with a certain complicated joint distribution for these
 limit quantities.

 To explain (\ref{torus-con}), first note that as $N \to \infty$ the
 differences
 $T^N_i - T^N_*$ are stochastically bounded (by the time to percolate
 through a finite set of edges) but cannot converge to $0$ (by linearity of
 growth rate in the shape theorem below), so we expect some
 non-degenerate
 limit distribution $(\tau_i , 1 \leq i \leq 4)$.
 Next consider the time $T^N_0$ at which the origin is wetted.
 By uniformity of starting position,
 $Q^N(T^N_0)$ must have uniform distribution on
 $\{1,2,\ldots,N^2\}$, and it follows that
 $N^{-2} Q^N(T^N_*) \cd U$.
 The final assertion
 \begin{equation}
  (N^{-1}(Q^N(T^N_* + t) - Q^N(T^N_*)), 0 \leq t < \infty)
 \cd
 (V t, 0 \leq t < \infty)
 \label{torus-con-2}
 \end{equation}
 is related to the {\em shape theorem} \cite{MR2045986}  for first-pasage
 percolation on the infinite lattice started at the origin.
 This says that the random set $\BB_s$ of vertices wetted before time $s$
 grows linearly with $s$, and the spatially rescaled set $s^{-1} \BB_s$ converges to
 a limit deterministic convex set $\BB$:
 \begin{equation}
 s^{-1} \BB_s \to \BB .
 \label{shape-theorem}
 \end{equation}
 It follows that
 \[ N^{-2} Q^N(sN) \to q(s) \mbox{ as } N \to \infty \]
 where $q(s)$ is the area of $s \BB$ regarded as a subset of the
 continuous
 torus $[0,1]^2$.
 Because
 $N^{-2} Q^N(T^N_0) \cd U$ we have
 \[ T^N_* \approx T^N_0 \approx N^2 q^{-1}(U) \]
 where $q^{-1}(\cdot)$ is the inverse function of $q(\cdot)$.
 Writing $Q^{\prime N}(\cdot)$ for a suitably-interpreted local growth
 rate
 of $Q^N(\cdot)$ we deduce
 \[ (N^{-2} Q^N(T^N_*), N^{-1}Q^{\prime N}(T^N_*)) \cd ( U,
 q^\prime(q^{-1}(U))) \]
 and so (\ref{torus-con-2}) holds for $V = q^\prime(q^{-1}(U))$.

 \subsection{Analysis of the rank-based reward game}
\label{sec-NN-ana}
 We want to study the case where other agents call some neighbor at rate
 $\theta$ but
 \ego\ (at the origin) calls some neighbor at rate $\phi$.
 To analyze rewards, by scaling time we can reduce to the case
  where other agents call {\em each} neighbor at rate $1$ and \ego\ calls each
 neighbor at rate $\lambda = \phi/\theta$.
 We want to compare the rank
 $M^N_\lambda$ of \ego\ (rank = $j$ if \ego\ is the $j$'th person to receive
 the information)
 with the rank $M^N_1$ of \ego\ in the $\lambda = 1$ case.
 As noted above, $M^N_1$ is uniform on $\{1,2,\ldots, N^2\}$.
 Writing $(\xi^\lambda_i, \ 1 \leq i \leq 4)$ for independent
 Exponential($\lambda)$ r.v.'s,
 the time at which the origin receives the information is
 \[ T^N_* + \min_i (T^N_i - T^N_* + \xi^\lambda_i) \]
 and the rank of the origin is
 \[ M^N_\lambda = Q^N(T^N_*) + N \widetilde{Q}^N(\min_i (T^N_i - T^N_* +
 \xi^\lambda_i)) \]
 where
 \[ \widetilde{Q}^N(t) = N^{-1}(Q^N(T^N_* + t) - Q^N(T^N_*)) . \] Note we
 can construct $(\xi^\lambda_i, \ 1 \leq i \leq 4)$ as
 $(\lambda^{-1} \xi^1_i, \ 1 \leq i \leq 4)$.
 Now use (\ref{torus-con-2}) to see that as $N \to \infty$
 \begin{equation}
  (N^{-2}M^N_1, N^{-1}(M^N_\lambda - M^N_1)) \cd
 (U, V Z(\lambda))
 \label{UVZ}
 \end{equation}
 where
\begin{equation}
  Z(\lambda) := \min_i (\tau_i + \xi^\lambda_i) - \min_i (\tau_i +
 \xi^1_i) . \label{Z-lam}
\end{equation}
 Now in the setting where \ego\ calls at rate $\phi$ and others at rate
 $\theta$
 we have
 \[
  \payoff(\phi,\theta) - \payoff(\theta,\theta) +
 (\phi - \theta)
 = E \left[ R \left(\frac{M^N_{\phi/\theta}}{N^2} \right)
 - R \left(\frac{M^N_{1}}{N^2} \right) \right]
 \]
 and it is straightforward to use (\ref{UVZ}) to show this
 \begin{equation}
  \sim N^{-1} \int_0^1 (-r(u)) \ z_u(\phi /\theta)  du
  , \mbox{ for } z_u(\lambda) := E(V Z(\lambda) | U=u).
 \label{payoff-NN}
 \end{equation}
 The Nash equilibrium condition
 \[
 \left. \frac{d}{d \phi} \payoff(\phi,\theta) \right|_{\phi = \theta} = 0
  \]
 now implies
 \begin{equation}
   \theta^{\mbox{{\tiny Nash}}}_N \sim N^{-1}
 \int_0^1 (-r(u)) \ z_u^\prime(1)  du .
 \label{zgu}
 \end{equation}
 Because $Z(\lambda)$ is decreasing in $\lambda$ we have  $z_u^\prime(1) <
 0$ and this expression is of the form (\ref{torus-result}) with 
 \begin{equation}
 g(u) = - z_u^\prime(1)  = - \sfrac{d}{d\lambda}  \left.E(VZ(\lambda)|U=u) \right|_{\lambda = 1}
\label{gz-formula}
\end{equation}

 \paragraph{Remark.}  The distribution of $V$ depends on the function
 $q(\cdot)$ which depends on the limit shape in nearest neighbor first passage percolation,
 which is not explicitly known.  Also
 $Z(\lambda)$ involves the joint distribution of $(\tau_i)$, which is not
 explicitly known, and also is
 (presumably) correlated with  the direction from the percolation source
 which is in turn not independent of $V$.  This suggests it would be
 difficult to find an explicit formula for $g(u)$.

\section{Order of magnitude arguments}
\label{sec-heuristics}
Here we mention simple order of magnitude arguments
for the two basic cases we have already analyzed.
As mentioned in the introduction, 
what matters is the size of the
 {\em window width} $w_{\theta,n}$ of the associated
  first passage percolation process
We will re-use such arguments in sections \ref{sec-grid-cN} and \ref{sec-5-heur}, in more complicated settings.  

\paragraph{Complete graph.}
If agents call at rate $\theta = 1$ then by (\ref{KS-lim}) the window width 
is order $1$; so if $\theta_n$ is the Nash equilibrium 
rate then the window width $w_n$ is order $1/\theta_n$.
Suppose $w_n \to \infty$.
Then \ego\ could call at some fixed slow rate $\phi$ and 
(because this implies many calls are made near the start of the window)
the reward to \ego\ will tend to $R(0)$, and \ego's payoff 
$R(0) - \phi$ will be larger than the typical payoff 
$\bar{R} - \theta_n$. 
This contradicts the definition of Nash equilibrium.
So in fact we must have $w_n$ bounded above, implying $\theta_n$ bounded below, implying the
Nash equilbrium in wasteful.

\paragraph{Nearest neighbor torus.}
If agents call at rate $\theta = 1$ 
then by the shape theorem (\ref{shape-theorem}) 
the window width is order $N$.
The time difference between receipt time for different neighbors of \ego\ 
is order $1$, so if \ego\ calls at rate $2$ instead of rate $1$ 
his rank (and hence his reward) increases by order $1/N$. 
By scaling, if the Nash equilibrium rate is $\theta_N$ 
and \ego\ calls at rate $2\theta_N$ then his increased reward is again of  order $1/N$.
His increased cost is $\theta_N$.
At the Nash equilibrium the increased reward and cost must balance, so 
$\theta_N$ is order $1/N$, so the Nash equilibrium is efficient.

  \section{The $N \times N$ torus with general interactions: a simple criterion for efficiency}
 \label{sec-grid-cN}

  {\bf Network communication  model.}
  The agents are at the vertices of the $N \times N$ torus.
 Each agent $i$ may, at any time, call any other agent $j$,
 at cost $c(d(i,j))$,
 and learn all items that $j$ knows.

 Here $d(i,j)$ is the distance between $i$ and $j$, and we assume the cost function
 $c(d)$ satisfies
 \begin{equation}
 c(1) = 1; \quad c(d) \uparrow \infty \mbox{ as } d \to \infty .
 \label{cd}
 \end{equation}

\noindent
 {\bf Poisson strategy.}
 An agent's strategy is described by a sequence
 $(\theta(d); \ d = 1,2,3,\ldots)$;  and for each $d$:\\
 \hspace*{0.7in}
 at rate $\theta(d)$ the agent calls a random agent at distance $d$.

\noindent
A simple argument below shows
 \begin{equation}
 \mbox{
 under condition (\ref{cd}) the Nash equilibrium is efficient.
 }
 \label{eff}
 \end{equation}
Consider the Nash strategy, and suppose first that the window width
 $w_N$ converges to a limit $w_\infty < \infty$.
 Consider a distance $d$ such that the Nash strategy has
 $\theta^{\mbox{{\tiny Nash}}}(d) > 0$.
 Suppose \ego\ uses $\theta(d) = \theta^{\mbox{{\tiny Nash}}}(d)  +
 \phi$.
 The increased cost is $\phi c(d)$ while the increased benefit is at most
 $O(w_\infty \phi)$, 
 because this is the increased chance of getting information earlier. 
 So the Nash strategy must have $\theta^{\mbox{{\tiny Nash}}}(d) = 0$ for
 sufficiently large $d$, not depending on $N$.
 But for first passage percolation with bounded range transitions,  the shape theorem (\ref{shape-theorem}) remains true and 
 implies that
 $w_N$ scales as $N$.

This contradiction implies that the window width $w_N \to \infty$.
 Now suppose the Nash equilibrium were inefficient, with some Nash cost
 $\bar{\theta} > 0$.
 Suppose \ego\ adopts the strategy of just calling a random neighbor at
 rate $\phi_N$,
 where $\phi_N \to 0, \ \phi_N w_N \to \infty$.
 Then \ego\ obtains asymptotically the same reward
 $\bar{R}$ as his neighbor, a typical agent.
 But \ego's cost is $\phi_N \to 0$.
 This is a  contradiction with the assumption of inefficiency.
 So the conclusion is that the Nash equilibrium is efficient and $w_N \to
 \infty$.
 
 \paragraph{Remarks.} 
 Result (\ref{eff}) is striking. but does not tell us what the Nash equilibrium strategy and cost actually are.  
 It is a natural open problem to study the case of (\ref{cd}) with $c(d) = d^\alpha$. 
 Instead we study a simpler model in the next section.

 \section{The $N \times N$ torus with short and long range interactions}
 \label{sec-torus-general}

  {\bf Network communication  model.}
  The agents are at the vertices of the $N \times N$ torus.
 Each agent $i$ may, at any time, call any of the $4$ neighboring agents
 $j$ (at cost $1$),
 or call any other agent $j$ at cost $c_N \geq 1$,
 and learn all items that $j$ knows.

\noindent
 {\bf Poisson strategy.}
 An agent's strategy is described by a pair of numbers
 $(\thnear, \thfar) = \theta$: \\
 \hspace*{0.7in}
 at rate $\thnear$ the agent calls a random neighbor \\
 \hspace*{0.7in}
 at rate $\thfar$ the agent calls a random non-neighbor.

 This model obviously interpolates between the complete graph model $(c_N
 = 1$) and the nearest-neighbor model ($c_N = \infty$).

 First let us consider for which values of $c_N$ the nearest-neighbor
 Nash
 equilibrium
 ($\thnear$ is order $N^{-1}$, \ $\thfar = 0$) persists in the current
 setting.
 When \ego\ considers using a non-zero value of $\thfar$, the cost is
 order
 $c_N \thfar$.
 The time for information to reach a typical vertex is order
 $N/\thnear = N^2$, and so the benefit of using a non-zero value of
 $\thfar$ is order
 $\thfar N^2$.
 We deduce that

 if $c_N \gg N^2$ then the Nash equilibrium is asymptotically the same as
 in the nearest-neighbor case;
 in particular, the Nash equilibrum is efficient.

 Let us study the more interesting  case
 \[ 1 \ll c_N \ll N^2 . \]
 The result in this case turns out to be, qualitatively 
 \begin{equation}
 \mbox{
 $\theta_{\mbox{{\tiny near}}}^{\mbox{{\tiny Nash}}}$ is order
 $c_N^{-1/2}$
 and
 $\theta_{\mbox{{\tiny far}}}^{\mbox{{\tiny Nash}}}$
 is order $c_N^{-2}$.
 In particular, the Nash equilibrum is efficient.}
 \label{torus-qual}
 \end{equation}
 ``Efficient" because the cost  $c_N \thfar + \thnear$ is order
 $c_N^{-1/2}$.  
 See (\ref{eq-99}) for the exact result.

We first do the order-of-magnitude calculation (section \ref{sec-5-heur}), then 
analyze the relevant first passage percolation process (section \ref{sec-5-fpp}), 
and finally do the exact analysis in section \ref{sec-5-analysis}.

 \subsection{Order of magnitude calculation}
 \label{sec-5-heur}
 Our order of magnitude argument for (\ref{torus-qual}) uses three ingredients
 (\ref{qual-1},\ref{qual-2},\ref{qual-3}).
As in section \ref{sec-heuristics}
we consider the
 {\em window width} $w_N$ of the associated percolation process. 
 Suppose \ego\ deviates from the Nash equilibrium
 $(\theta_{\mbox{{\tiny near}}}^{\mbox{{\tiny Nash}}},
 \theta_{\mbox{{\tiny
 far}}}^{\mbox{{\tiny Nash}}})$
 by setting his $\theta_{\mbox{{\tiny far}}} = \theta_{\mbox{{\tiny
 far}}}^{\mbox{{\tiny Nash}}} + \delta$.
The chance of thereby learning the information earlier,
and hence
the increased reward to \ego, is order $\delta w_N$
 and the increased cost is $\delta c_N$.
 At the Nash equilibrium these must balance, so
 \begin{equation}
 w_N \asymp c_N  \label{qual-1}
 \end{equation}
where $\asymp$ denotes ``same order of magnitude".
 Now consider the difference $\ell_N$ between the times that different
 neighbors of \ego\ are wetted.
 Then $\ell_N$ is order $1/\theta_{\mbox{{\tiny near}}}^{\mbox{{\tiny
 Nash}}}$.
 Write $\delta = \theta_{\mbox{{\tiny near}}}^{\mbox{{\tiny Nash}}}$ and
 suppose \ego\ deviates from the Nash equilibrium
 by setting his $\theta_{\mbox{{\tiny near}}} = 2 \delta$.
 The increased benefit to \ego\ is order $ \ell_N /w_N$
 and the increased cost is $\delta $.
 At the Nash equilibrium these must balance, so $\delta \asymp  \ell_N
 /w_N$ which becomes
 \begin{equation}
 \theta_{\mbox{{\tiny near}}}^{\mbox{{\tiny Nash}}}  \asymp w_N^{-1/2}
 \asymp c_N^{-1/2} . \label{qual-2}
 \end{equation}
 Finally we need to calculate how the window width $w_N$ for FPP depends
 on
 $(\thnear, \thfar)$, and we show in the next section that
 \begin{equation}
 w_N \asymp \thnear^{-2/3} \thfar^{-1/3} .
  \label{qual-3}
 \end{equation}
 Granted this, we substitute (\ref{qual-1},\ref{qual-2}) to get
 \[ c_N \asymp c_N^{1/3} \thfar^{-1/3} \]
 which identifies $\thfar \asymp c_N^{-2}$ as stated at
 (\ref{torus-qual}).

\subsection{First passage percolation on the $N \times N$ torus with short and long range interactions}
\label{sec-5-fpp}
We study the model  
(call it {\em short-long FPP}, to distinguish it from nearest-neighbor FPP)
defined by rates
\begin{eqnarray*}
\nu_{ij} &=& \sfrac{1}{4},\quad \quad  \ j \mbox{ a neighbor of } i \\
&=& \lambda_N/N^2, \ j \mbox{ not a neighbor of } i
\end{eqnarray*}
where $1 \gg \lambda_N \gg N^{-3}$.

Recall the shape theorem (\ref{shape-theorem}) for nearest neighbor first passage percolation; 
let $A$ be the area of the limit shape $\BB$.
Define an artificial distance $\rho$ such that $\BB$ is the unit ball in $\rho$-distance; so 
nearest neighbor first passage percolation moves at asymptotic speed $1$ with respect to 
$\rho$-distance.
Consider short-long FPP started at a random vertex of the $N \times N$ torus.
Write $F_{N,\lambda_N}$ for the proportion of vertices reached by time $t$ and let 
$T_{(0,0)}$ be the time at which the origin is reached.  
The event
$\{T_{(0,0)} \leq t\}$
corresponds asymptotically
to the event that at some time $t - u$ there is percolation 
across some long edge $(i,j)$ 
into some vertex $j$ at $\rho$-distance $\leq u$ from $(0,0)$   
(here we use the fact that 
nearest neighbor first passage percolation moves at asymptotic speed $1$ with respect to 
$\rho$-distance).
The rate of such events at time $t - u$ is approximately 
\[ N^2 F_{N,\lambda_N}(t-u) \times Au^2 \times \lambda_N/N^2   \]
where the three terms represent
the number of possible vertices $i$,
the number of possible vertices $j$,
and the percolation rate $\nu_{ij}$.
Since these events occur asymptotically as a Poisson process in time, we get
\begin{equation}
1 - F_{N,\lambda_N}(t) 
\approx P(T_{(0,0)} \leq t) 
\approx \exp \left(-A \lambda_N \int_0^\infty u^2 F_{N,\lambda_N}(t-u) \ du \right) . 
\label{T00N}
\end{equation}
This motivates study of the equation (for an unknown distribution function $F_\lambda$)
\begin{equation}
1 - F_\lambda(t) = \exp \left( - \lambda \int_{- \infty}^t (t-s)^2 F_\lambda(s) \ ds \right),
\quad - \infty < t < \infty 
\label{F-quad}
\end{equation}
whose solution should be unique up to centering.
Writing $F_1$ for the $\lambda = 1$ solution,
the general solution scales as 
\[
F_\lambda(t) := F_1(\lambda^{1/3}t) .
\]
So by (\ref{T00N}),
up to centering
\begin{equation}
 F_{N,\lambda_N}(t) \approx F_1((A \lambda_N)^{1/3}t)  . 
 \label{3scale}
 \end{equation}
To translate this result into the context 
of the rank-based
rewards game, 
suppose each agent uses strategy 
$\theta_N = (\thNnear,\thNfar)$.
Then the spread of one item of information is as 
first passage percolation with rates
\begin{eqnarray*}
\nu_{ij} &=& \thNnear /4, \ j \mbox{ a neighbor of } i \\
&=& \thNfar /(N^2-5), \ j \mbox{ not a neighbor of } i .
\end{eqnarray*}
This is essentially the case above with $\lambda_N = \thNfar/\thNnear$,
time-scaled by $\thNnear$, and so by (\ref{3scale})
the distribution function 
$F_{N,\theta_N}$ for the time at which a typical agent receives the information is
\begin{equation}
F_{N,\theta_N}(t) \approx F_1 \left( A^{1/3} \thNfar^{1/3} \thNnear^{2/3} t \right) . 
\label{F1-scale}
\end{equation}
In particular the window width is as stated at (\ref{qual-3}).

\subsection{Exact equations for the Nash equilibrium}
\label{sec-5-analysis}
The equations will involve three quantities:
\\
(i) The solution $F_1$ of (\ref{F-quad}).
\\
(ii) The area $A$ of the limit set $\BB$ in the 
shape theorem (\ref{shape-theorem}) for nearest-neighbor first passage pecolation.
\\
(iii) The limit distribution 
(cf. (\ref{torus-con}))
\begin{equation}
(T^r_i - T^r_*, \ 1 \leq i \leq 4)
\cd
(\tau_i, \ 1 \leq i \leq 4) 
\mbox{ as } r \to \infty
\label{tau-r}
\end{equation}
for relative receipt times of neighbors of the origin in 
nearest-neighbor first passage pecolation,
where now we start the percolation at a random vertex of $\rho$-distance 
$\approx r$ from the origin.

To start the analysis, 
suppose all agents use rates
$\theta = (\thNnear, \thNfar)$.
Consider the quantities

$S$ is the first time that \ego\ receives the information from a non-neighbor

$T$ is the first time that \ego\ receives the information from a neighbor

$F = F_{N,\theta_N}$ is the distribution function of 
$T$.

\noindent
With probability $\to 1$ as $N \to \infty$ \ego\ will actually receive the information first from a neighbor, and so $F$ is asymptotically 
the distribution function of the time at which \ego\ receives the information.

Now suppose \ego\ uses a different rate 
$\phiNfar \neq \thNfar$ 
for calling a non-neighbor.  
This does not affect $T$ but changes the distribution of $S$ to 
\[ P(S > t) \approx \exp\left( - \phiNfar \int_{-\infty}^t F(s) \ ds \right)
\]
by the natural Poisson process approximation.
Because $\thNfar$ is small we can approximate
\[ P(S \leq t) \approx \phiNfar \int_{-\infty}^t F(s) \ ds . \]
The mean reward to \ego\ for one item, as a function of $\phiNfar$, varies 
as
\[
E(R(F(S)) - R(F(T))1_{(S<T)} \ + \mbox{ constant.}
  \]
Because $U = F(T)$ is uniform on $(0,1)$, in the $N \to \infty$ limit
\begin{eqnarray*}
  E(R(F(S)) - R(F(T))1_{(S<T)}
&=& 
E(R(F(S)) - R(U))1_{(F(S) < U)}\\
&=& \int_0^1 \ du \ E(R(F(S)) - R(u))1_{(F(S) < u)} \\
&=& \int_0^1 \ du \ E \int_{\min(F(S),u}^u r(y) dy \\
&=& \int_0^1 dy \ (1-y) r(y) P(F(S) \leq y) \\
&=& \int_0^1 dy \ (1-y) r(y) P(S \leq F^{-1}(y)) \\
&=& \phiNfar \int_0^1 dy \ (1-y) r(y) \int_{-\infty}^{F^{-1}(y)} F(s) ds .
\end{eqnarray*}
The cost associated with using $\phiNfar$ is 
$c_N \phiNfar$, 
and at the Nash equilibrium the cost and reward must balance, 
so at the Nash equilibrium 
$F = F_{N,\theta_N}$ must satisfy
\begin{equation}
c_N \sim \int_0^1 dy \ (1-y) r(y) \int_{-\infty}^{F^{-1}(y)} F(s) ds .
\label{eq-1}
\end{equation}

Now suppose instead that \ego\ uses a different rate 
$\phiNnear \neq \thNnear$ 
for calling a neighbor.  
As in section \ref{sec-NN-ana},
we set $\lambda = \phiNnear/\thNnear$ 
so that we can use rate-$1$ 
nearest-neighbor first passage pecolation 
as comparsion.
For $(\tau_i)$ at (\ref{tau-r})
and independent Exponential($\lambda$) 
random variables
$(\xi^\lambda_i)$ write 
(as at (\ref{Z-lam}))
 \[ Z(\lambda) := \min_i (\tau_i + \xi^\lambda_i) - \min_i (\tau_i +
 \xi^1_i) . \]
So $Z(\lambda)$ is the time difference for \ego\ receiving the information, 
caused by \ego\  using $\phiNnear$ instead of $\thNnear$.
This time difference is measured after time-rescaling; in real time units 
the time difference is 
$Z(\lambda) / \thNnear$.

As above, 
write $T$ for receipt time for \ego\ using $\thNnear$, and $F = F_{N,\theta_N}$ for its distribution function.
Then receipt time for \ego\ using
$\phiNnear$
is
$T + Z(\lambda) / \thNnear$,
so \ego's rank becomes 
$\approx F(T) + F^\prime(T) Z(\lambda) / \thNnear$,
and setting $U = F(T)$ 
the rank of \ego\ is
$\approx U + F^\prime(F^{-1}(U)) Z(\lambda) / \thNnear$.
The associated mean reward change for \ego\  is asymptotically
\[ \sfrac{z(\lambda)}{\thNnear} \times 
\int_0^1 r(u) F^\prime(F^{-1}(u)) \ du ; \quad 
\lambda = \phiNnear/\thNnear  \]
where $z(\lambda) = EZ(\lambda)$.
Because the cost of using rate
$\phiNnear$ 
equals 
$\phiNnear$,
the Nash equilibrium condition (\ref{Nash-crit}) implies
\begin{equation}
\thNnear^2 \sim z^\prime(1) \int_0^1 r(u) F^\prime(F^{-1}(u)) \ du  .
\label{eq-2}
\end{equation}

We have now obtained the desired two equations for 
$F_{N,\theta_N}$ at the Nash equilibrium $\theta_N$.
Use (\ref{F1-scale}) to rewrite these equations (\ref{eq-1},\ref{eq-2}) in terms of $F_1$ as 
\begin{eqnarray*}
c_N &\sim& 
A^{-1/3} \thNfar^{-1/3} \thNnear^{-2/3} 
\int_0^1 dy \ (1-y) r(y) \int_{-\infty}^{F_1^{-1}(y)} F_1(s) ds \\
\thNnear^2 &\sim&  A^{1/3} \thNfar^{1/3} \thNnear^{2/3}
z^\prime(1) 
\int_0^1 r(u) F_1^\prime(F_1^{-1}(u)) \ du .
\end{eqnarray*}
Solving for 
$\thNnear, \thNfar$ we find 
\begin{equation}
\thNnear \sim Q^{1/2} c_N^{1/2}, \quad 
\thNfar \sim A^{-1} Q^{-1} c_N^{-2} 
\label{eq-99}
\end{equation}
for
\[ Q = 
z^\prime(1) 
\left(\int_0^1 dy \ (1-y) r(y) \int_{-\infty}^{F_1^{-1}(y)} F_1(s) ds \right)
\left(\int_0^1 r(u) F_1^\prime(F_1^{-1}(u)) \ du \right)
 . \]

 \section{Variants}
 \subsection{Transitivity and the symmetric variant}
 \label{sec-transitive}
 The examples we have studied so far have a certain property called {\em
 transitivity} in graph theory \cite{big74}.
 Informally, transitivity means ``the network looks the same to each
 agent";
 formally, it means that for any two agents $i,j$ there is an
 automorphism
 of the network that preserves the network cost structure and maps $i$ to
 $j$.
 This is what allows us to assume that  in a Nash equilibrium each agent
 uses same strategy.

 The general framework of section  \ref{sec-frame} uses the
  {\em asymmetric} model in which
 agent $i$ calls  agent $j$ (at a certain cost to $i$)  and learns all
 items that $j$ knows.
 In the {\em symmetric} variant, agent $i$ calls  agent $j$ (at a certain cost
 to $i$),
 and each tells the other all items they know.

 For the transitive networks we have studied there is a simple
 relationship
 between the Nash equilibrium values of the asymmetric and symmetric
 variants of the Poisson strategies:
 \begin{equation}
  \theta^{\mbox{{\tiny Nash}}}_{\mbox{{\tiny sym}}}  = \sfrac{1}{2}
   \theta^{\mbox{{\tiny Nash}}} _{\mbox{{\tiny asy}}} .
   \label{asy-sym}
   \end{equation}
   The point is that the percolation process in the symmetric variant is
 just the
   percolation process in the asymmetric variant, run at twice the speed,
 and this leads to
   the
   following relationship between the reward when \ego\ uses rate $\phi$ and other agents use rate $\theta$:
   \[ \reward_{\mbox{\tiny sym}}(\phi,\theta) = \reward_{\mbox{\tiny
 asy}}(\phi + \theta, 2 \theta) . \]
   Because $\payoff(\phi,\theta) = \reward(\phi,\theta) - \phi$ in each
 case, we get
   \[ \payoff_{\mbox{\tiny sym}}(\phi,\theta) = \payoff_{\mbox{\tiny
 asy}}(\phi + \theta, 2 \theta) + \theta  \]
   and therefore
     \[ \frac{d}{d\phi} \payoff_{\mbox{\tiny sym}}(\phi,\theta) =
 \frac{d}{d\phi} \payoff_{\mbox{\tiny asy}}(\phi + \theta, 2 \theta) . \]
 The criterion (\ref{Nash-crit}) leads to (\ref{asy-sym}).
 
 \subsection{Communication at regular intervals}
 \label{sec-reg}
 We have studied  ``Poisson rate $\theta$" calling strategies because these are simplest to analyze explicitly.  
 A natural alternative is the 
  ``regular, rate $\theta$" strategy in which agent $i$ calls a random
 other
 agent at times
 \begin{equation}
 U_i, U_i + \sfrac{1}{\theta}, U_i + \sfrac{2}{\theta}, \ldots
 \end{equation}
 where $U_i$ is uniform on $(0,\frac{1}{\theta})$.  
 
 Consider first the complete graph case, and the setting (section \ref{sec-finite}) 
  where (for fixed $k \geq 2$) the first $k$
 recipients get reward
 $n/k$.
In this case, for $k = 2$ formula (\ref{k=2ego}) is replaced by
 \[ P( \mbox{\ego\  is second to receive item}) =
 \int_0^{\min(\frac{1}{\phi},\frac{1}{\theta})} (1 - \theta u)^{n-2} \
 \phi
 \ du  \]
 and repeating the analysis in section \ref{sec-finite} gives exactly the same asymptotics
 (\ref{nash-2-1},\ref{nash-2-2})
 as in the Poisson case.  
 Consider instead the general reward framework
 \[ \mbox{  The $j$'th person to learn an item of information gets reward
 $R(\frac{j}{n})$} .\]
 If all agents use rate $\theta$ then the distribution function $F_\theta$ for receipt time for a typical agent satisfies 
 (as an analog of the logistic equation (\ref{eq-logistic-theta})) 
 \begin{equation}
 1 - F_\theta(t) = \int_{t - \frac{1}{\theta}}^t \ \prod_{i \geq 0} \left (1-F_\theta(s - \sfrac{i}{\theta}) \right) \ \theta \ ds .
 \end{equation}
 If \ego\  switches to rate $\phi$ then the distribution function $G_{\phi,\theta}$ for \ego's receipt time 
 satisfies (as an analog of (\ref{Gphit}))
  \begin{equation}
 1 - G_{\phi,\theta}(t) = \int_{t - \frac{1}{\phi}}^t \ \prod_{i \geq 0} \left (1-F_\theta(s - \sfrac{i}{\phi}) \right) \ \phi \ ds .
 \end{equation}
 One can now continue the section \ref{sec-CG-analy} analysis; 
 we do not get useful explicit solutions but the qualitative behavior is similar to the 
 ``Poisson calls" case, and in particular the Nash equilibrium is wasteful.
 
Similarly, on the $N \times N$ grid with nearest neighbor interaction, switching from the ``Poisson calls"
case to the ``regular calls" case preserves the order $N^{-1}$ value of the Nash equilibrium rate 
$\theta^{\mbox{{\tiny Nash}}}_N$ and hence preserves its efficiency.

 \subsection{Gossip with reward based on audience size}
 \label{sec-audience}
Perhaps a more realistic model for gossip is to replace Rule 2 by

\noindent
{\bf Rule 3.} An agent $i$ gets reward $c$ whenever another agent learns an item from $i$.

\noindent 
For the complete graph and Poisson($\theta$) strategies we can re-use the section \ref{sec-CG-analy} analysis 
to calculate the Nash equilibrium.  
First suppose all agents use the same rate $\theta$ and consider an agent $i$ who receives the information at percentile $u$.  
For $j > un$ the $j$'th agent to receive the information has chance  $\frac{1}{j}$ to 
receive it from agent $i$, and so the mean reward to agent $i$ is 
(calculations in the $n \to \infty$ limit) 
$c \int_u^1 \frac{1}{x} \ dx = -c \log (1-u)$.  
Suppose now \ego\ switches to rate $\phi$.  
Then (calls incur unit cost) 
\[ \payoff(\phi,\theta) = - \phi + c E(- \log (1 - F_\theta(T_{\phi,\theta}))) \]
where the time $T_{\phi,\theta}$ at which \ego\ receives the information has distribution function $G_{\phi,\theta}$ at (\ref{eq-G}), and where $F_\theta$ at (\ref{F-theta}) is the distribution function of the time at which a typical agent receives the information.  
Now
\begin{eqnarray*}
E(- \log (1 - F_\theta(T_{\phi,\theta}))) 
&=& \int_0^1 \sfrac{1}{1-u} P(F_\theta(T_{\phi,\theta}) \leq u) \ du \\
&=& \int_0^1 \sfrac{1}{1-u} (1 - (1-u))^{\phi/\theta} \ du \mbox{ by } (\ref{FTu}) \\
&=& \int_0^1 y^{-1} (1 - y)^{\phi/\theta} \  dy
\end{eqnarray*}
and then we calculate
\[  \frac{d}{d \phi} \payoff(\phi,\theta) = - 1 - c \int_0^1 \sfrac{\log (1-y)}{\theta} \ \sfrac{(1-y)^{\phi/\theta}}{y} \  dy . \] 
Now the Nash equilibrium criterion (\ref{Nash-crit}) implies 
\begin{equation}
\theta^{\mbox{{\tiny Nash}}}_n \to  - c \int_0^1 \sfrac{1-y}{y} \ \log (1-y) \ dy .
\end{equation} 
So switching to this ``Rule 3" model preserves the wastefulness of the Nash equilibrium on the compete graph.

However, for the $N \times N$ grid with nearest neighbor interaction, 
switching to the ``Rule 3" models changes the efficient 
($\theta^{\mbox{{\tiny Nash}}}_N$ is order $N^{-1}$) Nash equilibrium to a wasteful equilibrium 
with $\theta^{\mbox{{\tiny Nash}}}_N$ becoming order $1$.

 \subsection{Related literature} 
 We do not know any literature closely related to our model.  
 As well as the epidemic and the gossip algorithm topics mentioned in 
the introduction, and classic applied probability work on {\em stochastic rumors} \cite{daley-kendall}, other loosely related work includes 
\begin{itemize}
\item 
models where agents form networks under
 conditions where there are costs for maintaining network edges and
 benefits from being part of a large network \cite{jackson2003smn}. 
 \item Prisoners' Dilemma games between neighboring agents on a graph \cite{hauert:405}.
 \end{itemize}
 One can add many other topics which are harder to model mathematically,
 e.g. diffusion of technological innovations  \cite{innovations} or of
 ideologies.


\end{document}